\documentclass{article}
\usepackage{graphicx} 
\usepackage{amsfonts}
\usepackage{amsthm}
\usepackage{amsmath}

\usepackage{authblk}
\usepackage{orcidlink}

\def\ff{\widehat f}

\def\gg{\widehat g}

\newcommand{\ee}{\mathbf{e}}

\def\erre {\mathbb{R}}

\def \be {  \varpi}

\newcommand{\fer}[1]{(\ref{#1})}

\def\bx{{\bf x}}
\def\bX{{\bf X}}
\def\bY{{\bf Y}}

\def\by{{\bf y}}

\def\bxi{{\pmb\xi}}

\newcommand{\R}{\mathbb R}
\def\be#1\ee{\begin{equation}#1\end{equation}}

\newenvironment{equations}{\equation\aligned}{\endaligned\endequation}

\newtheorem{theorem}{Theorem}
\newcommand{\rev}{\textcolor{black}}

\title{ On the domain of attraction of source--type solutions of the porous medium equations}

\author[1]{Gennaro AURICCHIO\orcidlink{0000-0002-4285-8887}}
\author[2]{Giuseppe TOSCANI\orcidlink{0000-0003-4419-804X}}

\affil[1]{Department of Mathematics, University of Padua, via Trieste 36, Padova, 35100 Italy}
\affil[2]{Department of Mathematics, University of Pavia, and IMATI of the National Council for Research, via Ferrata 1, Pavia, 27100 Italy}

\begin{document}

\maketitle

\begin{abstract}
    We study the rate of convergence to Barenblatt profiles of the solutions  to  porous medium  equations  in $\R^n, n \ge 1$,  by means of a certain class of Energy distances, a class of discrepancies between probability measures which have been recently widely used in problems related to AI, in particular for tasks related to machine learning. If  $n \ge 3$, it is shown both that the rate of decay  in Energy distance is improved with respect to known results which are based on the decay of relative entropy, and that, with respect to this distance,  any initial value of bounded mass belongs to the domain of attraction of the source--type solution. Similar results are shown to hold for the fast diffusion case.
\end{abstract}

\section{Introduction}\label{intro}

The study of the asymptotic behavior of the solutions $f(\bx,t)$ of the nonlinear diffusion equation 
\be\label{eq:por}
\frac{\partial f(\bx,t)}{\partial t} = \Delta f^m(\bx,t), \qquad m > \frac{n-2}n
\ee
where $\bx \in \R^n, n \ge 1$, is a classical problem which, starting with the porous medium case, characterized by $m>1$, has been originally studied  in the \rev{$70$'s} of the last century in a pioneering paper \rev{by} Kamin \cite{Kam} (cf. \cite{Vaz,Vaz1} for an exhaustive presentation of the  results obtained on this topic up to the recent years).
In a physical sense, the solution $f(\cdot, t)$ to \eqref{eq:por} represents the time-density  of the diffusion process.
If we consider the Cauchy Problem for equation \fer{eq:por} posed in $\R^n \times \R_+$ with integrable and nonnegative initial data 
\be\label{eq:ini}
f_0(\bx) \in L^1(\R^n), \quad f_0(\bx) \ge 0
\ee
the source-type solutions of equation \fer{eq:por} \rev{provide} the asymptotic representation of a large class of solutions starting from $f_0$ satisfying \fer{eq:ini} \cite{Vaz}.

The source-type or Barenblatt solutions of equation \fer{eq:por} are given by \cite{Bar}
\be\label{eq:Bar}
b(\bx,t, C) = t^{-n\alpha}F_C\left(t^{-\alpha}\bx \right)
\ee
where $C$ is a positive free constant, and $F_C$, the profile, has the form
\be\label{eq:det}
F_C(\by) = \left( C-\kappa |\by|^2\right)^{1/(m-1)}_+, \quad \alpha = \frac 1{n(m-1) +2}, \quad \kappa = \frac{m-1}{2m}\alpha. 
\ee
In \fer{eq:det}, the function   {$(g)_+$}
denotes as usual the positive part of $g$.

In this range of the parameter $m$, for a given initial value $f_0$ of given mass $M$, the  constant $C = C(f_0)$ can be chosen to satisfy the rule of mass conservation
\be\label{eq:mc}
M= \int_{\R^n} f_0(\bx)\, d\bx = \int_{\R^n} b(\bx,t, C)\, d\bx.
\ee
Further, in presence of additional conditions on the initial density, such as    \rev{boundedness of the relative (to the Barenblatt solution $F_C$) Ralston-type entropy of order $m$ \cite{CaTo1}, or the boundedness of the relative Renyi entropy of order $m$ \cite{CaTo,SaTo}},  it is possible to reckon an explicit rate of convergence to the right source-type solution in $L^1(\R^n)$. It is interesting to note that the boundedness of energy allows to establish a direct connection between the  asymptotic stabilization property of the porous medium flow and the stabilization to the Gaussian profile embodied in the classical central limit theorem \cite{To1}. Moreover, the boundedness of energy is essential to characterize in a precise way the best matching Barenblatt profile \cite{DoTo1,DoTo2} 

Diffusion processes described by \fer{eq:por}  arise in various fields, such as plasma physics, kinetic theory of gases, solid state physics, filtration models and transport in porous media. In some cases, including many metals and ceramic materials the diffusion coefficient can, over a wide range of temperatures, be approximated as $|\bx|^{-\nu}$ for some
exponent $\nu$ such that $0 < \nu < 2$ \cite{Ros}, which leads to an exponent $m < 1$ in \fer{eq:por}. Also, the special case $\nu = 1$ emerges  in the central limit approximation to Carleman’s model of the Boltzmann equation, according to \cite{LiTo,McK}. It is interesting to remark that, while less relevant from the point of view of applications, in mathematics the fast diffusion case $m<1$ is better known than the porous medium case for its properties in connection with functional analysis.

In what follows, we approach the problem of identifying an explicit rate of convergence of the solution to equation \fer{eq:por} to the source-type solution in a new way, which  requires only that the initial datum has a finite  mass. This broadened range of attraction of the Barenblatt solution can be obtained only if the distance in which the explicit rate of convergence is proven is weaker than the classical $L^1$-distance; moreover, while as usual the rate of convergence is closely dependent on the space dimension $n$,  the proposed approach works in this generality only in dimension $n \ge 3$.

To this end, we make use of arguments similar to the ones recently introduced in \cite{AuTo}, which are based on applying Cramér distance \cite{Cra} to one-dimensional Fokker–Planck equations with linear drift and nonlinear diffusion. As extensively discussed  in \cite{AuTo}, the Cramér distance—also known as the continuous ranked probability score—has gained renewed attention in AI-related problems, particularly in machine learning \cite{Belle}, where it is used for tasks such as distribution estimation, regression, and image generation.
As first remarked in \cite{Belle}, the main motivation for adopting the Cramér distance, a metric for comparing one-dimensional probability distributions, lies in its advantageous properties, most notably its unbiased sample gradients, which make it especially amenable to optimization via stochastic gradient descent.

As we shall detail in this paper, the results proven in \cite{AuTo} for one-dimensional Fokker--Planck type equations follow by resorting to properties typical of  Cramér's distance, closely linked to some of the various forms in which this distance can be expressed. As a matter of fact, an interesting related distance, which, as shown in \cite{ABGT2}, shares most properties of Cramér's distance, while designed to also compare multivariate probability distributions, is the Energy distance introduced by Székely \cite{Sze89, Sze03, SR} by resorting to \rev{one of the alternative expressions of Cram\'er’s distance}.  

The classical Energy distance, of wide use in Statistics,  has been recently generalized in \cite{ABGT2}, where the relationships with other well-known metrics—such as the Wasserstein distance \cite{villani2009optimal} and the Fourier-based metrics \cite{CaTo,GTW}—were thoroughly analyzed, in reason of their applicability in the field of artificial intelligence \cite{Au1,Au2,cuturi2014fast,scagliotti2025normalizing}, as well as in the context of partial differential equations \cite{For0,For1,For2}.

In this paper, we show that a certain class of the new Energy--type distances considered in \cite{ABGT2} are very-well adapted to study the rate of convergence to the source--type solutions to non linear diffusion equations. 
{These metrics are directly connected with the Cram\'er distance in the one-dimensional case.
However,} at variance with the one-dimensional Cram\'er distance, this class of $n$-dimensional Energy distances  do not work well on drift operators for $n >1$, and \rev{cannot} be directly applied to study the large-time behavior of multi-dimensional Fokker--Planck equations.

In more details, in Section \ref{sec:basic} we briefly report the definitions and the main properties of both Cram\'er and Energy distances, as well as their connections. The equivalence between these two distances in the one-dimensional case relies on  arguments first applied to  obtain \rev{an} alternative definition of Gini index \cite{Xu}.  Then, in Section \ref{sec:energy} we will show how a particular class of these distances can be fruitfully employed to study convergence to source--type solutions to the porous medium equation \fer{eq:por}, departing from any initial value of bounded mass.
Section \ref{sec:decay} will be devoted to \rev{proving} that an explicit rate of convergence follows as soon as $n \ge 3$.  
In Section \ref{sec:fast}, we \rev{will} show that a similar argument can be used to study the convergence rate to source--type solutions in the fast diffusion case $m<1$ if $m$ belongs to a certain range and if $n \ge 3$.
{In Section \ref{sec:conclusions}, we \rev{will} outline the conclusions of the paper.}

\section{From Cram\'er to Energy distances}\label{sec:basic}

Let $f$ be a probability density over $\R^n, n \ge 1$. 
We will associate $f$ \rev{with} a random variable $\textbf{X}$, such that for a subset of the reals $A \subseteq \R^n$, we have 
\[
P(\textbf{X}\in A) = \int_A f(\bx)\, d\bx.
\]
Given $s>0$, we denote by $P_s(\R^n, n \ge 1)$ the class of all probability densities $f$ over $\R^n$ such that
\[
m_s(f)  = \int_{\R^n} |\bx|^s f(\bx)\, d\bx < + \infty.
\]
If $s=1$,  the quantity $m_1(f)$ defines the (absolute) mean value of $\textbf{X}$, say $E(|\textbf{X}|)$. 

Given a probability density $f$ in $P_s(\R^n), n \ge 1$, its Fourier transform is defined by
\[
\ff(\bxi) = \int_{\R}  e^{-2\pi i\bxi\cdot\bx} \, f(\bx)\, d\bx.
\]
If $n =1$, the (cumulative) distribution function of \rev{the one-dimensional random variable $X$} is given by
\be\label{cumu}
F(x) = P(X \le x) = \int_{-\infty}^x f(x)\, dx.
\ee
 Given two \rev{one-dimensional} random variables $X$ and $Y$  {and their} cumulative distribution functions $F$  {and $G$}, the Cram\'er distance between $X$ and $Y$ is expressed by \cite{Cra}
\be\label{Cramer}
d(F,G) = \int_\R[F(x) -G(x)]^2\, dx.
\ee 
Note that, as written, Cram\'er distance is not a metric, while its square root is. 

One of the interesting properties of Cram\'er's distance is that it can be equivalently written  in different ways. It is immediate to show that it can be written in terms of Fourier transforms.
Indeed, owing to Parseval's theorem, Cram\'er distance \fer{Cramer} can be equivalently written as
 \be\label{Cramer2}
d(F,G)  =  \int_\R \frac{|\ff(\xi) -\gg(\xi)|^2}{|\xi|^2} \, d\xi.
\ee
Also, assuming that the \rev{one-dimensional} random variables $X$ and $Y$ are mutually independent, and that $X'$ and $Y'$ are {independent copies of $X$ and $Y$}, Cram\'er's distance can be expressed as
\be\label{Cramer3}
d(F,G) = E(|X-Y|) - \frac 12 E(|X-X'|) - \frac 12 E(|Y-Y'|).
\ee
Note that \fer{Cramer3} implies that the Cram\'er distance of $F$ and $G$ is bounded as soon as $F,G \in P_1(\R)$. This third way of expressing Cram\'er's distance has been shown in \cite{SR} (cf. \cite{AuTo} for a detailed proof based on the Gini index \cite{Gini1,Gini2}).

{Expressing the Cram\'er distance as in \fer{Cramer3} allows \rev{us} to compute a weak metric between two probability distributions in closed form on every $\R^n$, with $n >1$.}
 In a series of papers  \cite{Sze89, Sze03, SR}, Székely introduced the definition of Energy distance between two multidimensional random variables $\bX$ and $\bY$ with densities $f(\bx)$, and, respectively, $g(\bx)$, with $\bx \in \R^n$,  $n \ge 1$. The Energy distance in $\R^n, n \ge 1$, of order $\alpha$, with $0<\alpha<2$ is defined by
\begin{equations}\label{energy-a}
E_{n,\alpha} (f,g) &= 2\int_{\R^n\times\R^n}|\bx-\by|^\alpha f(\bx)g(\by) \, d\bx d\by\\
& -\int_{\R^n\times\R^n}|\bx-\by|^\alpha f(\bx)f(\by) \, d\bx d\by
 -\int_{\R^n\times\R^n}|\bx-\by|^\alpha g(\bx)g(\by) \, d\bx d\by\\
&= - \int_{\R^n\times\R^n}|\bx-\by|^\alpha [f(\bx)-g(\bx)][f(\by) -g(\by)] \, d\bx d\by.
\end{equations} 
As for the Cram\'er distance, \rev{ in dimension $n>1$} the Energy distance of order $\alpha$ can be equivalently expressed in Fourier transform. In this case it holds \cite{ABGT2,SR}
\be\label{eq:energy2}
E_{n,\alpha} (f,g) = c_{n,\alpha} \int_{\R^n} \frac{ |\ff(\bxi) -\gg(\bxi)|^2}{|\bxi|^{n+\alpha}}\, d\bxi,
\ee
{where the constant $c_{n,\alpha}$, determined as in \cite{Stein}} takes the value
\be\label{cn}
c_{n,\alpha} = \alpha 2^{n+\alpha-1} \pi^{n/2} \frac{\Gamma\left( \frac{n+\alpha}2\right)}{\Gamma\left( \frac{2-\alpha}2\right)}.
\ee
\rev{where $\Gamma$ is the Gamma function}.

The passage to Fourier transform is the easiest way to verify that right-hand side in \fer{energy-a} is nonnegative, so that {twice} the mixed integral dominates the sum of the last two. Further, expression \fer{energy2} allows one to relate in a precise way the energy distance of order $\alpha$ to homogeneous Sobolev spaces of fractional order. 

Let us note that in one dimension, and for $\alpha=1$, the Energy distance  $E_{1,1}(f,g)$ is nothing but twice the Cram\'er distance, as defined by \fer{Cramer3}, so that
\begin{equations}\label{eq:chiave}
d(f,g) &= \frac 12 E_{1,1}(f,g)
 = 2\int_{\R\times\R}\frac{|x-y|}2 f(x)g(y) \, dx dy \\
& -\int_{\R\times\R}\frac{|x-y|}2 f(x)f(y) \, dx dy
 -\int_{\R\times\R}\frac{|x-y|}2 g(x)g(y) \, dx dy = \\
 & -\int_{\R\times\R}\frac{|x-y|}2 [f(x)-g(x)][f(y)-g(y)] \, dx dy.
\end{equations}

 \section{Learning  from Cram\'er distance}\label{sec:energy}
In what follow, we will show how Cram\'er distance works in connection with the porous medium equation \fer{eq:por}.  We assume that the initial value for the Cauchy problem for \fer{eq:por} satisfies conditions \fer{eq:ini}, and all \rev{solutions satisfy} well-known regularity properties  that are enough to justify all the computation that follows \cite{Vaz0}.  

A key property of identity \fer{eq:chiave} is that  for any given probability density $h \in P_1(\R)$ it holds
 \be\label{eq:d2}
 \frac{d^2}{dx^2}  \int_\R \frac{|x-y|}2 h(y) \,  dy = h(x),
 \ee
 namely that the function 
 \[H(x)=  \frac{d^2}{dx^2}  \frac{|x|}2
 \]
  coincides with the delta function $\delta_0(x)$ centered in $x=0$. 
  
  This implies that, if $f(x,t)$ and $g(x,t)$ are the two probability density functions \rev{solving} the one-dimensional porous medium equation \fer{eq:por} departing from $f_0(x)$ (respectively $g_0(x)$) satisfying \fer{eq:ini} with mass equal to one, and such that $d(f_0,g_0)$ is bounded, integrating twice by parts one has
  \begin{equations}\label{chiave}
&\frac d{dt}d(f(t),g(t)) =\frac 12  \frac d{dt}E_{1,1}(f(t),g(t)) =\\
  &- \frac d{dt} \int_{\R\times\R}\frac{|x-y|}2 [f(x,t)-g(x,t)][f(y,t)-g(y,t)] \, dx dy=\\
  & - 2 \int_{\R\times\R}\frac{|x-y|}2 \frac d{dt} [f(x,t)-g(x,t)][f(y,t)-g(y,t)] \, dx dy = \\
  &  - 2 \int_{\R\times\R}\frac{|x-y|}2[f(y,t)-g(y,t)]\frac {d^2}{dx^2} [f^m(x,t)-g^m(x,t)] \, dx dy =\\
   &  - 2 \int_\R[f^m(x,t)-g^m(x,t)]\frac {d^2}{dx^2} \left\{\int_\R \frac{|x-y|}2[f(y,t)-g(y,t)] \, dy \right\}\, dx = \\
    &  - 2 \int_\R[f^m(x,t)-g^m(x,t)][f(x,t)-g(x,t)] \,   dx  \le 0.
\end{equations}
 
Hence, the Cram\'er distance between any two solutions to the porous medium equation is decreasing in time. 

The key point in deriving this result is that the kernel $|x|/2$ in the definition of Cram\'er distance satisfies \fer{eq:d2}.

This result can be easily generalized to any dimension $n >1$ by resorting to the concept of fundamental solution of the Laplacian $\Delta$ in $\R^n$, namely a function $\Phi(\bx)$ that satisfies in $\R^n$  the equation $\Delta \Phi(\bx) = \delta_0(\bx)$ in the sense of distributions \cite{Smo}. The fundamental solutions read
\be\label{eq:fund}
\Phi_n(x)=
\begin{cases}
\displaystyle \frac{1}{(2-n)\,|S_{n-1}|}\,|\bx|^{2-n}, & n\ge 3,
\\
\displaystyle \frac{1}{2\pi}\,\log{|\bx|}, & n=2,\\
\end{cases}
\qquad \bx\in\mathbb{R}^n\setminus\{0\},
\ee
where $|S_{n-1}|$ denotes the surface area of the unit sphere in $\R^n$, given by
\[
|S_{n-1}|=\frac{2\pi^{n/2}}{\Gamma\!\left(\tfrac{n}{2}\right)}.
\]
Owing to \fer{eq:fund} we introduce, for $n \ge 3$,  a (normalized) variant of the Energy distance between two $n$-dimensional probability densities  in the form
\begin{equations}\label{energy-n}
&\mathcal E_{n} (f,g) = 2\int_{\R^n\times\R^n}\frac{|\bx-\by|^{2-n}}{(2-n)\,|S_{n-1}|} f(\bx)g(\by) \, d\bx d\by\\
& -\int_{\R^n\times\R^n}\frac{|\bx-\by|^{2-n}}{(2-n)\,|S_{n-1}|} f(\bx)f(\by) \, d\bx d\by
 -\int_{\R^n\times\R^n}\frac{|\bx-\by|^{2-n}}{(2-n)\,|S_{n-1}|} g(\bx)g(\by) \, d\bx d\by\\
&{= \int_{\R^n\times\R^n}\frac{|\bx-\by|^{2-n}}{(n-2)\,|S_{n-1}|} [f(\bx)-g(\bx)][f(\by) -g(\by)] \, d\bx d\by.}
\end{equations} 
The Energy distance of negative order $2-n$, with $n\ge 3$, defined in \fer{energy-n} is a particular case of the energy distances with negative exponent introduced and studied in \cite{ABGT2}.
{As shown in in \cite{ABGT2}, the distance $\mathcal E_n$ is, up to a multiplicative constant $c_{n,\alpha}$ given in \fer{cn}, equal to a Fourier \rev{based} metric, which in this case gives} 

\be\label{energy2}
\mathcal E_{n} (f,g) =  \int_{\R^n} \frac{ |\ff(\bxi) -\gg(\bxi)|^2}{|\bxi|^{2}}\, d\bxi,
\ee
{so that $\mathcal E_{n}(f,g)$ is the square of the distance in the Sobolev space $\dot H_{-1}(\R^n)$.}

\section{The decay in Energy distance}\label{sec:decay}

Let us now consider two probability density functions  $f(\bx,t)$ and $g(\bx,t)$, solutions of the $n$-dimensional porous medium equation \fer{eq:por} departing from $f_0(\bx)$ (respectively $g_0(\bx)$)  and such that $\mathcal E_n(f_0,g_0)$ is bounded. By means of Green's third identity, we can repeat the computations in \fer{eq:chiave} for any dimension $n>1$, obtaining the decay
\be\label{eq:chiave-n}
 \frac d{dt}\mathcal E_n(f(t),g(t)) =  - 4 \mathcal I_n(f(t),g(t)),
\ee
where
\be\label{eq:ep}
\mathcal I_n(f,g) =  \int_{\R^n}[f^m(\bx)-g^m(\bx)][f(\bx)-g(\bx)] \, d\bx \ge 0.
\ee
Similar computations have been done in \cite{AuTo} in the (simpler) linear case $m=1$ for the classical Energy distance corresponding to $0<\alpha<2$.  

 {In the next step, we} obtain a lower bound on the functional $\mathcal I_n$ in terms of the energy distance $\mathcal E_n$. 

First of all, since $m>1$,   inequality 
\be\label{vero}
x^m-a^m  -(x-a)^m \ge 0, \quad x\ge a \ge 0
\ee
implies
\be\label{eq:ep1_0}
\mathcal I_n(f,g) \ge  \int_{\R^n}|f(\bx)-g(\bx)|^{1+m} \, d\bx.
\ee
Further, since $m >1$, thanks to the classical interpolation inequality
\be\label{eq:inter}
\int_{\R^n} |h(\bx)|^2 \, dx \le \left[\int_{\R^n} |h(\bx)| \, dx\right]^{(m-1)/m}\left[\int_{\R^n} |h(\bx)|^{1+m} \, dx\right]^{1/m},
\ee
and to the fact that $f$ and $g$ are two probability densities on $\R^n$, so that
\[
 \int_{\R^n}|f(\bx)-g(\bx)|\, d\bx \le \int_{\R^n}[f(\bx)+g(\bx)]\, d\bx = 2,
\]
we get the inequality
\be\label{eq:ine-m}
\mathcal I_n(f,g)  \ge  \frac 1{2^{m-1}}\left[ \int_{\R^n}|f(\bx)-g(\bx)|^2 \, d\bx \right]^m,
\ee
which shows that, as $m>1$,  the functional $\mathcal I_n(f,g)$ controls the $m$-th power of the square of the $L^2$-distance between $f$ and $g$.

Now, consider that, for $n \ge 3$,  we can bound from above the Energy distance  $\mathcal E_n(f,g)$ between two probability densities by their $L^2$-distance resorting to a Nash-type argument.  For any given positive constant $R$ we  have
\begin{equations}\label{eq:Nash}
& \int_{\R^n} \frac{ |\ff(\bxi) -\gg(\bxi)|^2}{|\bxi|^{2}}\, d\bxi = \\& \int_{|\bxi| \le R} \frac{ |\ff(\bxi) -\gg(\bxi)|^2}{|\bxi|^{2}}\, d\bxi + \int_{|\bxi|>R} \frac{ |\ff(\bxi) -\gg(\bxi)|^2}{|\bxi|^{2}}\, d\bxi  \le\\
&  \int_{|\bxi| \le R} \frac{ 2^2}{|\bxi|^{2}}\, d\bxi + \frac 1{R^2}  \int_{\R^n} |\ff(\bxi) -\gg(\bxi)|^2\, d\bxi  = \\
& {4|S_{n-1}|}\int_{\rho\le R} \frac{\rho^{n-1}}{\rho^{2}}\, d\rho + \frac 1{R^2}  \int_{\R^n} |f(\bx) -g(\bx)|^2\, d\bx=\\
& C_n \, R^{n-2}  + \frac 1{R^2}  \int_{\R^n} |f(\bx) -g(\bx)|^2\, d\bx,
\end{equations}
where
\be\label{eq:cn}
C_n = \frac{8 \pi^{n/2}}{(n-2)\Gamma\left( \frac n2\right)}
\ee
Optimizing with respect to  $R\in \R_+$ we obtain the bound
\be\label{eq:Nash-o}
\mathcal E_n(f,g) \le d_n \left[\int_{\R^n} |f(\bx) -g(\bx)|^2\, d\bx\right]^{(n-2)/n},
\ee
where the constant $d_n$ is explicit
\be\label{eq:dn}
d_n= \frac n{n-2}\frac{2^{4/n}\pi}{\Gamma\left( \frac n2\right)^{2/n}}  .
\ee
Owing to \fer{eq:ine-m} we finally have
\be\label{eq:Nash-I}
\mathcal E_n(f,g) \le d_n 2^{(m-1)/m}\mathcal I_n(f,g)^{(n-2)/nm}.
\ee
Substituting into \fer{eq:chiave-n} we obtain that the time derivative of the Energy distance $\mathcal E_n$ satisfies the inequality
\be\label{eq:chiave-fin}
 \frac d{dt}\mathcal E_n(f(t),g(t)) \le   - D(m,n) \mathcal E_n(f(t),g(t)) ^{\frac{mn}{n-2}},
\ee
where the explicit constant $D(m,n)$ is given by
\[
D(m,n) = \frac 4{d_n 2^{(m-1)m}}.
\]
The differential inequality \fer{eq:chiave-fin} can be integrated to give the time-decay
\be\label{eq:En-decay}
\mathcal E_n(f(t),g(t)) \le \left[\left(\frac 1{\mathcal E_n(f_0,g_0)}\right)^{\gamma(m,n)}+ D(m,n)\gamma(m,n)\,t \right]^{-1/\gamma(m,n)}.
\ee
In \fer{eq:En-decay}
\be\label{eq:delta}
\gamma(m,n) = \frac{(m-1)n +2}{n-2}.
\ee
Note that, for any given $m >1$ the value of $\gamma(m,n)$ is decreasing with $n\ge 3$, and
\be\label{eq:b-a}
m-1 \le \gamma(m,n) \le 3m-1.
\ee
being the lower value $m-1>0$ reached as $n \to +\infty$.

Hence we proved
\begin{theorem}\label{the:decay}
Let $n \ge 3$, and let  $f(\bx,t)$ and $g(\bx,t)$  be the solutions of the $n$-dimensional porous medium equation \fer{eq:por} departing from the initial values $f_0(\bx)$ (respectively $g_0(\bx)$) of equal mass $M=1$ such that the Energy distance $\mathcal E_n(f_0,g_0)$ is bounded. Then, $\mathcal E_n(f(t),g(t))$ decays to zero as $t \to +\infty$ at the explicit rate $(1/t)^{1/\gamma(m,n)}$, where $\gamma(m,n)$ is defined in \fer{eq:delta}. 
\end{theorem}
 
 {We now apply Theorem \ref{the:decay} to the case in which the initial configuration is the Dirac delta measure centered in zero}.
Since the solution to the porous medium equation departing from \rev{$\delta_0$} is nothing but the Barenblatt solution {$b(t)$} of mass one defined in \fer{eq:Bar}, we get both the convergence to zero and the rate of decay to the source--type solution in the Energy distance $\mathcal E_n(f(t),b(t))$, provided the Energy distance is bounded at the initial time $t=0$. 
 {Indeed, when this condition is met,} $f_0(\bx)$ belongs to the domain of attraction of the source--type solution $b(\bx,t,C)$ with respect to the Energy distance $\mathcal E_n$. 

Applying definition \fer{energy-n} to the pair $(f_0,\delta_0)$  {yields}
\begin{equations}\label{energy-0}
0 \le \mathcal E_{n} (f_0,\delta_0) = 2&\int_{\R^n}\frac{|\bx|^{2-n}}{(2-n)\,|S_{n-1}|} f_0(\bx) \, d\bx \\
& -\int_{\R^n\times\R^n}\frac{|\bx-\by|^{2-n}}{(2-n)\,|S_{n-1}|} f_0(\bx)f_0(\by) \, d\bx d\by\\
=&{\int_{\R^n\times\R^n}\frac{|\bx-\by|^{2-n}}{(n-2)\,|S_{n-1}|} f_0(\bx)f_0(\by) \, d\bx d\by} \\
& -2\int_{\R^n}\frac{|\bx|^{2-n}}{(n-2)\,|S_{n-1}|} f_0(\bx) \, d\bx.
\end{equations} 
Since the Energy distance is non-negative, and $n \ge 3$, $\mathcal E_{n} (f_0,\delta_0)$ is bounded as soon as
\[
\int_{\R^n\times\R^n}{|\bx-\by|^{2-n}} f_0(\bx)f_0(\by) \, d\bx d\by
\]
 is bounded. 
 
Owing to the Hardy--Littlewood--Sobolev inequality \cite{Lieb}, stating that
 \[
\left| \int_{\R^n\times\R^n}{|\bx-\by|^{-\lambda}} h(\bx)k(\by) \, d\bx d\by\right| \le N_{p,\lambda,n} \|h\|_p\|k\|_r
 \]
for all $h \in L^p(\R^n)$, $k \in L^r(\R^n)$ with  $1 < p,r  < \infty$, $1/p + 1/r + \lambda/n = 2$  and $0 < \lambda < n$, an explicit sufficient condition  \rev{ensuring} the boundedness of the Energy distance in \fer{energy-0} is to assume that $f_0\in L^p(\R^n)$ with $p = 2n/(n+2)$.

Let the initial density $f_0(\bx) \in L^1(\R^n)$, and let $\ff_0(\bxi)$ its Fourier transform. Also,
let $\omega(\bx)$ denote the Gaussian density
\be\label{eq:Gau}
\omega(\bx) = \frac 1{\sqrt{2\pi}}\exp\left\{ -\frac{|\bx|^2}2\right\}.
\ee
Owing to expression \fer{energy2}, and using the bound $|\ff_0(\bxi)| \le 1$, we obtain
\begin{equations}\label{eq:L1-bound}
\mathcal E_{n} (f_0,\delta_0) &=  \int_{\R^n} \frac{ |\ff_0(\bxi) -1|^2}{|\bxi|^{2}}\, d\bxi = \\
& \int_{\R^n} \frac{ \left|\ff_0(\bxi) -\ff_0(\bxi)\widehat\omega(\bxi) + \ff_0(\bxi)\widehat\omega(\bxi)  -1\right|^2}{|\bxi|^{2}}\, d\bxi \le \\
& 2 \int_{\R^n} \frac{\left|\widehat\omega(\bxi)-1 \right|^2|\ff_0(\bxi)|^2}{|\bxi|^{2}}\, d\bxi + 
2  \int_{\R^n} \frac{ \left|\ff_0(\bxi)\widehat\omega(\bxi)  -1\right|^2}{|\bxi|^{2}}\, d\bxi \le \\
&2 \int_{\R^n} \frac{\left|\widehat\omega(\bxi)-1 \right|^2}{|\bxi|^{2}}\, d\bxi + 
2 \int_{\R^n} \frac{ \left|\ff_0(\bxi)\widehat\omega(\bxi)  -1\right|^2}{|\bxi|^{2}}\, d\bxi \le C <+\infty.
\end{equations}
Indeed, since both $\omega$ and $f_0\star\omega$, the convolution of $f_0$ and $\omega$, belong to $L^1(\R^n) \cap L^\infty(\R^n)$,  they also belong to $L^p(\R^n)$ when $p = 2n/(n+2)$.
{In particular,} their Energy distances from the Dirac delta located in zero, according to   the Hardy--Littlewood--Sobolev inequality, are bounded.

Therefore we have
\begin{theorem}\label{the:decay-to-B}
Let $n \ge 3$, and let  $f(\bx,t)$ be the solution of the $n$-dimensional porous medium equation \fer{eq:por} departing from a probability density $f_0(\bx)$.  Then, the Energy distance $\mathcal E_n(f(t),b(t,C))$, {with $b(t,C)$ \rev{defined} as in \eqref{eq:Bar},} decays to zero as $t \to +\infty$ at the explicit rate $(1/t)^{1/\gamma(m,n)}$, where $\gamma(m,n)$ is defined in \fer{eq:delta}. 
\end{theorem}

Theorem \ref{the:decay-to-B} shows that,  { if $n\ge 3$}, the domain of attraction of the source--type solution { defined as in \fer{eq:Bar} with respect to the Energy distance $\mathcal E_n$} contains all initial values in  $L^1(\R^n)$.
Hence, the domain of attraction of the Barenblatt solution is larger than the one found in \cite{CaTo0} by resorting to entropy methods, since there convergence in the stronger $L^1$-distance requires the initial density to have the second moment bounded. 
{Moreover, we stress that in \cite{CaTo0} the convergence of the solution to the porous medium equation towards the source--type solution in $L^1(\R^n)$ was shown with rate $(1/t)^{1/\mu(m,n)}$, where} 
\[
\mu(m,n) = \frac{(n+2)m-n}2. 
\]
{so that,  for $n =1,2$, we have $\gamma(n,m) < \mu(n,m)$.
However, if $n=3$, the two rates coincide, while if $n>3$, the rate of convergence in Energy distance is better than the rate of convergence in $L^1$, and  does not collapse to zero if $n \to\infty$.}

\section{The fast diffusion case}\label{sec:fast}

Let us now consider the fast diffusion case. To avoid confusion with the previous Sections, let us denote with $p<1$ the characteristic exponent of the nonlinear diffusion \fer{eq:por} in the fast diffusion regime. Moreover, since we want to study the large time behavior of the solution to equation \fer{eq:por}  {in which} the Barenblatt profile has bounded mass, we assume that $p >(n-2)/n$. It is immediate to realize that all the computations of Section \ref{sec:energy} leading to  equation \fer{eq:chiave-n}, with $\mathcal I_n(f,g)$ defined in \fer{eq:ep} continue to hold. Hence, as in Section \ref{sec:decay} we still obtain that the Energy distance $\mathcal E_n$ satisfies
\be\label{eq:chiave-nf}
 \frac d{dt}\mathcal E_n(f(t),g(t)) =  - 4 \mathcal J_n(f(t),g(t)),
\ee
where now, for $\frac{n-2}{n}<p<1$, 
\be\label{eq:fd}
\mathcal J_n(f,g) =  \int_{\R^n}[f^p(\bx)-g^p(\bx)][f(\bx)-g(\bx)] \, d\bx \ge 0.
\ee
In this case, however, we can not proceed as in Section \ref{sec:decay}, since inequality \fer{vero} is true if and only if $m>1$. On the other hand, for any given $m >1$, inequality \fer{vero} can be rewritten in the form
\[
x-a  -(x^{1/m}-a^{1/m})^m \ge 0, \quad x\ge a \ge 0.
\]
namely, since $1/m =p <1$ 
\be\label{vero1}
x-a  -(x^p-a^p)^{1/p} \ge 0, \quad x\ge a \ge 0
\ee
that implies
\be\label{eq:ep1}
\mathcal J_n(f,g) \ge  \int_{\R^n}|f^p(\bx)-g^p(\bx)|^{1+1/p} \, d\bx.
\ee
We now provide a lower bound on the right-hand side of \eqref{eq:ep1}.

Let $f$ and $g$  two probability density functions and $p\in[1/2,1)$. Thanks to the pointwise inequality
\[
|f(\bx)-g(\bx)|\le   |f^p(\bx)-g^p(\bx)||f^{1-p}(\bx)+g^{1-p}(\bx)|,
\]
for any given positive constant $\beta$, \rev{the} Holder inequality allows \rev{us}  to conclude that
\begin{equation}
\label{eq:betainequality}
    \int_{\erre^n}|f(\bx)-g(\bx)|^\beta d\bx\le\Big(\int_{\erre^n}|f^p(\bx)-g^p(\bx)|^{\beta r}d\bx\Big)^{\frac{1}{r}}\Big(\int_{\erre^n}|f^{1-p}(\bx)+g^{1-p}(\bx)|^{\beta s}d\bx\Big)^{\frac{1}{s}},
\end{equation}
with $r$ and $s$ such that ${1}/{r}+{1}/{s}=1$.
We then \rev{look} for $\beta, r$, and $s$ such that
\be\label{eq:valori}
    \frac{1}{r}+\frac{1}{s}=1 ; \quad
    \beta r= \frac{1+p}{p}; \quad
    \beta s(1-p)=1.
\ee
Solving the system, we find the values
\be\label{eq:beta}
\beta=\frac{1+p}{1+p-p^2}; \quad  r=\frac{1+p-p^2}{p}; \quad s=\frac{1+p-p^2}{1-p^2}.
\ee
Note that  $r, s>1$, while $1<\beta<2$.

Since $\beta s(1-p) =1$, and $f,g$ are probability densities, owing to the
classical inequality
\[
(a^q+b^q)^{1/q} \le 2^{1+1/q} (a +b),
\]
valid for any given $0<q<1$ and $a,b>0$, we obtain
\[
\int_{\erre^n}|f^{1-p}(\bx)+g^{1-p}(\bx)|^{\beta s}d\bx \le  2^{1+1/(1-p)}\int_{\erre^n}(f(\bx)+g(\bx))\, d\bx = 2^{(3-p)/(1-p)}.
\]
Substituting into \fer{eq:betainequality}, and taking into account \fer{eq:ep1} we finally get, for any $p \in (1/2,1)$
\be\label{eq:key-beta}
\mathcal J_n(f,g) \ge c_p\left(\int_{\erre^n}|f(\bx)-g(\bx)|^\beta d\bx\right)^r,
\ee
where $\beta=\beta(p)$ and $r=r(p)$ are given in \fer{eq:valori}, and
\be\label{eq:cp}
c_p= 2^{-(3+2p-p^2)/p}.
\ee
Hence $\mathcal J(f,g)$ controls the $L^\beta$-distance of the probability densities $f$ and $g$. 

Since $\beta$ in \fer{eq:beta} satisfies $1<\beta <2$,  the  Beckner form of \rev{the} Haussdorf--Young inequality \cite{Bec}
\[
 \left(\int_{\erre^n}|\ff(\bxi)-\gg(\bxi)|^{q^*} d\bxi\right)^{1/q^*}
\le \left(\frac{{q}^{1/q}}{{q^*}^{1/q^*}}\right)^{n/2}\left(\int_{\erre^n}|f(\bx)-g(\bx)|^q d\bx\right)^{1/q},
\]
which holds true for $1<q \le 2$ and $1/q+1/q^*=1$, allows \rev{us} to bound 
the  $L^{\beta^*}$-distance of the Fourier transforms of the probability densities $f$ and $g$. 

Since $\beta r = (1+p)/p$, \rev{the} Haussdorf--Young inequality implies 
\be\label{eq:HY}
\left(\int_{\erre^n}|f(\bx)-g(\bx)|^\beta d\bx\right)^r\ge \left(\frac{{\beta^*}^{1/\beta^*}}{\beta^{1/\beta}}\right)^{n(1+p)/(2p)}\left(\int_{\erre^n}|\ff(\bxi)-\gg(\bxi)|^{(1+p)/p^2} d\bxi\right)^p.
\ee
Therefore, inequality \fer{eq:key-beta} allows \rev{us} to conclude with the following 
\be\label{eq:key-p}
\mathcal J_n(f,g) \ge c_{p,n} \left(\int_{\erre^n}|\ff(\bxi)-\gg(\bxi)|^{(1+p)/p^2} d\bxi\right)^p,
\ee
where, clearly, the exponent $(1+p)/p^2 > 2$, while
\be\label{eq:cpn}
c_{p,n} = c_p \left(\frac{{\beta^*}^{1/\beta^*}}{\beta^{1/\beta}}\right)^{n(1+p)/(2p)}.
\ee
Thanks to \fer{eq:key-p}, we can proceed as in Section \ref{sec:decay} to obtain
\begin{align}
    \int_{\erre^n}\frac{|\hat f (\bxi) -\hat g(\bxi)|^2}{|\bxi|^2}d\bxi&\le \int_{|\bxi|\le R}\frac{|\hat f(\bxi)-\hat g(\bxi)|^2}{|\bxi|^2}d\bxi+\int_{|\bxi|\ge R}\frac{|\hat f(\bxi)-\hat g(\bxi)|^2}{|\bxi|^2}d\bxi \\
    &\le C_n R^{n-2}+\int_{|\bxi|\ge R}\frac{|\hat f(\bxi)-\hat g(\bxi)|^2}{|\bxi|^2}d\bxi,
\end{align}
where $R>0$ is a given value and $C_n$ is the constant  in \fer{eq:cn}.

By \rev{the} Holder inequality, we can estimate the right hand side of the latter inequality as it follows
\begin{align}
\label{eq:ineq_1}
    \int_{|\bxi|\ge R}\frac{|\hat f(\bxi)-\hat g(\bxi)|^2}{|\bxi|^2}d\bxi \le \bigg(\int_{|\bxi|\ge R}\frac{1}{|\bxi|^{2u}}d\bxi\bigg)^{{1}/{u}}\bigg(\int_{\erre^n}|\hat f(\bxi)-\hat g(\bxi)|^{2v}d\bxi\bigg)^{1/v},
\end{align}
where $u$ and $v$ are such that ${1}/{u}+{1}/{v}=1$.
{
If $v=\frac{1+p}{2p^2}>1$, then $u=\frac{1+p}{1+p-2p^2}$, so that, making use of \fer{eq:key-p} we obtain
\begin{align}
\label{eq:ineq_2}
    \int_{|\bxi|\ge R}\frac{|\hat f(\bxi)-\hat g(\bxi)|^2}{|\bxi|^2}d\bxi \le \bigg(\int_{|\bxi|\ge R}\frac{1}{|\bxi|^{\frac{2(1+p)}{1+p-2p^2}}}d\bxi\bigg)^{\frac{1+p-2p^2}{1+p}}\left(\frac 1{c_{p,n}} \mathcal J_n(f,g)\right)^{\frac{2p}{1+p}}.
\end{align}
Notice that $|\bxi|^{-\frac{2(1+p)}{1+p-2p^2}}$ is integrable on the set $|\bxi| \ge R$ if and only if 
\[
\frac{2(1+p)}{1+p-2p^2} -n + 1 > 1, 
\]
or, equivalently,
\be\label{eq:inep}
n(2p^2-p-1) + 2(1+p) >0.
\ee
It is immediate to verify that, to solve \fer{eq:inep} on $\R_+$,  $p$ must be greater than 
\[
p(n) = \frac{n-2}{4n}+\frac{1}{4n}\sqrt{9n^2-20n+4}.
\]
Notice that if $n=3$, then $p(3)={1}/{2}$.
Moreover, if $n\to\infty$, then $p(n)\to 1$.
For such values of $p$, we have
\begin{equation}
    \bigg(\int_{|\bxi|\ge R}\frac{1}{|\bxi|^{\frac{2(1+p)}{1+p-2p^2}}}d\bxi\bigg)^{\frac{1+p-2p^2}{1+p}}\le D_{p,n}R^{\frac{n(1+p-2p^2)}{1+p}-2}=D_{p,n}R^{-\nu(p,n)},
\end{equation}
where we denoted 
\be\label{eq:nupn}
\nu(p,n):=2- \frac{n(1+p-2p^2)}{1+p}>0,
\ee
and
\be\label{eq:dpn}
D_{p,n} =\left( \frac{|S_{n-1}|}{\frac{2(1+p)}{1+p-2p^2} -n}\right)^{\frac{1+p-2p^2}{(1+p)}}.
\ee
}

\rev{All in all}, we get
\begin{align}
\label{eq:ineq_3}
    \int_{\erre^n}\frac{|\hat f (\bxi) -\hat g(\bxi)|^2}{|\bxi|^2}d\bxi \le C_n R^{n-2}+ D_{p,n}c_{p,n}^{-\frac{2p}{1+p}} R^{-\nu(p,n)}\mathcal J_n(f,g)^{\frac{2p}{1+p}}.
\end{align}
By optimizing the right hand side of \eqref{eq:ineq_3} with respect to $R$ we get
\begin{align}
\label{eq:ineq_3b}
    \int_{\erre^n}\frac{|\hat f (\bxi) -\hat g(\bxi)|^2}{|\bxi|^2}d\bxi \le K(p,n)\mathcal{J}_n(f,g)^{\frac{n-2}{np}},
\end{align}
where the constant $K(p,n)$, that can be explicitly computed in terms of $C_n$, $c_{p,n}$, $\nu(p,n)$ and  $D_{p,n}$, as given by  \fer{eq:cn},  \fer{eq:cpn},  \fer{eq:nupn} and \fer{eq:dpn}, takes the value
\be\label{eq:kpn}
K(p,n) = \frac{2np^2}{1+p} \left(\frac{(n-2)C_n}{\nu(p,n)}\right)^{\frac{\nu(p,n)(1+p)}{2np^2}} D_{p,n}^{\frac{(n-2)(1+p)}{2np^2}} c_{p,n}^{-\frac{n-2}{np}}
\ee
What is important to remark is that, apart from the different value of the constant $K(p,n)$, the previous computations show that in the fast diffusion case the time decay of the energy distance $\mathcal{E}_n(f,g)$ follows the law \fer{eq:chiave-fin} relative to the porous medium equation. 

If we now go back to the notation $p=m$, where now  $m<1$, we conclude that the decay is given by  \fer{eq:En-decay}, where now $\gamma(m,n)$ is well defined only when the fast diffusion coefficient $m <1$ satisfies the bound
\be\label{eq:well}
m > \max\left\{\frac 12, \frac{n-2}n\right\}.
\ee
Note that in this case, the measure of the interval $((n-2)/n, 1)$ in which the coefficient $\gamma(m,n)$ is well defined converges to zero as $n \to \infty$.

\rev{All in all}, we proved the following.

\begin{theorem}\label{the:decayfast}
Let $n \ge 3$, and let  $f(\bx,t)$ and $g(\bx,t)$  be the solutions of the $n$-dimensional fast diffusion equation \fer{eq:por} departing from the initial values $f_0(\bx)$ (respectively $g_0(\bx)$) of equal mass $M=1$ such that the Energy distance $\mathcal E_n(f_0,g_0)$ is bounded. Then, if $m$ satisfies condition \fer{eq:well}, $\mathcal E_n(f(t),g(t))$ decays to zero as $t \to +\infty$ at the explicit rate {$(1/t)^{1/\gamma(m,n)}$, where $\gamma(m,n)$ is defined in \fer{eq:delta}.} 
\end{theorem}

Likewise, since the reasoning leading to Theorem \ref{the:decay-to-B}, and in particular the estimate \fer{eq:L1-bound}, do not depend on the exponent of the diffusion equation \fer{eq:por}, but only on the properties of the initial value, we have the following.

\begin{theorem}\label{the:decay-to-Bfast}
Let $n \ge 3$, and let  $f(\bx,t)$   be the solution of the $n$-dimensional fast diffusion equation \fer{eq:por} departing from a probability density $f_0(\bx)$.  Then, if $m$ satisfies condition \fer{eq:well}, the Energy distance $\mathcal E_n(f(t),b(t,C))$ decays to zero as $t \to +\infty$ at the explicit rate {$(1/t)^{1/\gamma(m,n)}$, where $\gamma(m,n)$ is defined in \fer{eq:delta}.} 
\end{theorem}

\section{Conclusions}\label{sec:conclusions}

In this paper,  we studied the rate of convergence  of the solution to  the porous medium equation  to the Barenblatt source--type solution with the same mass  by means of a particular class of Energy distances, which are equivalent to the square of the distance in the Sobolev space $\dot H_{-1}(\R^n)$.
It is shown that, as $n \ge 3$, the rate of decay is improved with respect to known results which are based on the decay of relative entropy, and that the domain of attraction of the Barenblatt solution includes all initial value of finite mass. 

With few modifications, it has been further shown that analogous results  hold in the fast diffusion case, with some limitations in the diffusion exponent that are mainly due to guarantee that the source--type solution has a bounded mass. 

The present results suggest that both Cram\'er and Energy distances could be of high interest in the study of both theoretical and numerical problems related to nonlinear diffusion equations.

\vfill\eject
\section{Acknowledgement}
This work has been written within the activities of GNFM (Gruppo Nazionale per la Fisica Matematica) of INdAM (Istituto Nazionale di Alta Matematica), Italy. G.T. acknowledges support of the Institute for Applied
Mathematics and Information Technologies (IMATI), Pavia, Italy. G.T. wishes to thank the Editors  for inviting him to contribute to this volume dedicated to Giuseppe Savar\'e, a friend and colleague at the Department of Mathematics at the University of Pavia for many years.


\begin{thebibliography}{99}


\bibitem{ABGT2} Auricchio, G., Brigati, G., Giudici, P., and  Toscani, G.,  From kinetic theory to AI: a rediscovery of high-dimensional divergences and their properties. \emph{Mathematical Models and Methods in Applied Sciences}, \textbf{36} (2026) 1--49.

\bibitem{Au1}
Auricchio, G., Codegoni, A., Gualandi, S., Toscani, G., and Veneroni, M.,
On the equivalence between Fourier-based and Wasserstein metrics, \textit{Rendiconti Lincei. Matematica e Applicazioni}, \textbf{31}  (2020) 627--649.

\bibitem{Au2}
Auricchio, G., Codegoni, A., Gualandi, S., and Zambon, L.,
The Fourier discrepancy function, \textit{Communications in Mathematical Sciences}, \textbf{21}  (2023) 627--639.

\bibitem{AuTo}
Auricchio, G., and  Toscani, G., Energy distance and evolution problems: a promising tool for kinetic equations. \emph{Ricerche di Matematica} \textbf{75} (2026) 849--863. 

\bibitem{Bar}
Barenblatt, G.I., On some unsteady motions of a liquid and gas in a porous medium. \emph{Akad. Nauk SSSR
Prikl. Mat. Meh.} \textbf{16} (1952) 67--78 

\bibitem{Bec}
Beckner, W.,  Inequalities in Fourier analysis. \emph{Annals of Mathematics}, \textbf{102} (1) (1975) 159--182 

\bibitem{Belle}
  Bellemare, M. G., Danihelka, I., Dabney, W., Mohamed, S., Lakshminarayanan, B., and Hoyer, S., and Munos, R., The Cram\'er distance as a solution to biased Wasserstein gradients, arXiv preprint arXiv:1705.10743 (2017).

  

\bibitem{CaTo0}
Carrillo, J.A. and Toscani, G. Asymptotic $L^1$-decay of solutions of the porous medium equation to self-similarity. \emph{Indiana University Mathematics Journal}, (2000) Apr 1:113--42.


\bibitem{CaTo}
Carrillo, J.A. and Toscani, G., Contractive probability metrics and asymptotic behavior of dissipative kinetic equations, \emph{Rivista di Matematica dell’Università di Parma},
\textbf{6} (2007) 75--198.

\bibitem{CaTo1}
Carrillo, J.A. and Toscani, G., Renyi entropy and improved equilibration rates to self-similarity for nonlinear diffusion equations, \emph{Nonlinearity}, \textbf{27}  (2014) 3159--3177.




\bibitem{Cra}
Cram\'er, H., On the composition of elementary errors, \emph{Skandinavisk Aktuarietidskrift}, \textbf{11}, (1928) 141--180.

\bibitem{cuturi2014fast}
Cuturi, M. and Doucet, A.,  Fast computation of Wasserstein barycenters, \emph{International conference on machine learning}, (2014) 685--693.

\bibitem{For0}
Di Francesco, M., Fornasier, M., Hu\"utter, J.C. and Matthes, D., Asymptotic behavior of gradient flows driven by nonlocal power repulsion and attraction potentials in one dimension. \emph{SIAM Journal on Mathematical Analysis}, \textbf{46} (6), (2014) 3814--3837.

\bibitem{DoTo1}
Dolbeault, J., and Toscani, G., Best matching Barenblatt profiles are delayed, \emph{J. Phys. A: Math. Theor.}, \textbf{48} (2015) 065206.

\bibitem{DoTo2}
Dolbeault, J., and Toscani, G., Nonlinear diffusions: extremal properties of Barenblatt profiles, best matching and delays. \emph{NonlinearAnal. Series A}, \textbf{138} (2016) 31--43.  


\bibitem{For1}
Fornasier, M., Haskovec, J. and Steidl, G.,  Consistency of variational continuous-domain quantization via kinetic theory. \emph{Applicable Analysis}, \textbf{92} (6), (2013)  1283--1298.

\bibitem{For2}
Fornasier, M. and Hütter, J.C.,  Consistency of probability measure quantization by means of power repulsion–attraction potentials. \emph{Journal of Fourier Analysis and Applications}, \textbf{22} (3) (2016) 694--749.



\bibitem{GTW}
Gabetta, E., Toscani, G., and Wennberg, B., Metrics for probability measures and the trend to equilibrium for solutions of the Boltzmann equation, \textit{Journal of Statistical Physics}, \textbf{81}  (1995) 901--934.


\bibitem{Gini1}
Gini, C.,  Sulla misura della concentrazione e della variabilit\`a dei caratteri, \textit{Atti del Reale Istituto Veneto di Scienze, Lettere ed Arti}, \textbf{73}  (1914) 1203--1248. English translation in \textit{Metron}, 3--38, (2005).

\bibitem{Gini2}
Gini, C.,  Measurement of inequality of incomes, \textit{The Economic Journal}, \textbf{31}  (1921) 124--126.

\bibitem{Kam}
Kamenomostskaya (Kamin), S., The asymptotic behaviour of the solution of the filtration equation,
\emph{Israel J. Math.} \textbf{14} (1973) 76–-87.



\bibitem{Lieb}
Lieb, E.H., Sharp constants in the Hardy-Littlewood-Sobolev and related inequalities, \emph{Ann. Math.} \textbf{118} (1983) 349--374.

\bibitem{LiTo}
Lions, P.L. and Toscani, G.,  Diffusive limit for finite velocity Boltzmann kinetic models. \emph{Revista Matematica Iberoamericana} \textbf{13} (3) (1997) 473--513.


  




\bibitem{McK}
McKean, H.P., The central limit theorem for Carleman’s equation. \emph{Israel Journal of Mathematics} \textbf{21} (1) (1975) 54--92.
 
\bibitem{Ros}
Rosenau, P.,  Fast and superfast diffusion processes. \emph{Physical review letters}, \textbf{74} (7)(1995) 1056--1059


\bibitem{SaTo}  
Savar\'e G., and  Toscani, G., The concavity of R\'enyi entropy power. \emph{IEEE Transactions on Information theory}, \textbf{60} (5) (2014) 2687--93.

\bibitem{scagliotti2025normalizing}
Scagliotti, A. and Farinelli, S., {Normalizing flows as approximations of optimal transport maps via linear-control neural ODEs}, \textit{Nonlinear Analysis},  \textbf{257} (2025), 113811.

\bibitem{Smo}
Smoller, J., \emph{Shock Waves and Reaction--Diffusion Equations} (2 ed.). Springer New York 1974. 

\bibitem{Stein}
Stein, E.M, and Weiss, G., \emph{Introduction to Fourier analysis on Euclidean spaces}, Princeton Univerdsity Press, Princeton, NJ (1970)

\bibitem{Sze89}
Sz\'ekely, G.J., \emph{Potential and kinetic energy in statistics}. Lecture Notes. Budapest Institute of Technology (Technical University) (1989).

\bibitem{Sze03}
Sz\'ekely, G.J.,  E-statistics: Energy of statistical samples. Bowling Green State University, Department of Mathematics and Statistics Technical Report No. 03-05, 2003 (also technical reports by same title, from 2000–2003 and NSA Grant \# MDA 904-02-1-0091 (2000–2002).)

\bibitem{SR}
Sz\'ekely, G.J. and Rizzo, M.I.,
Energy statistics: A class of statistics based on distances. \emph{Journal of Statistical  Planning and Inference}, \textbf{143} (2013) 1249--1272




\bibitem{To1}
Toscani, G., A central limit theorem for solutions of the porous medium equation, \emph{J. Evol. Equ.} \textbf{5} (2005) 185--203 

\bibitem{Vaz}
V\'azquez, J.L, Asymptotic behaviour for the porous medium equation posed in the whole space. \emph{Nonlinear Evolution Equations and Related Topics: Dedicated to Philippe Bénilan}. Basel: Birkhäuser Basel, 2003. 67--118.

\bibitem{Vaz0}
V\'azquez, J.L., \emph{The porous medium equation: mathematical theory}. Clarendon Press, 2006.

\bibitem{Vaz1}
V\'azquez, J.L., The mathematical theories of diffusion: nonlinear and fractional diffusion. In: \emph{Nonlocal and Nonlinear Diffusions and Interactions: New Methods and Directions: Cetraro, Italy 2016}. Cham: Springer International Publishing, 2017. 205--278.



\bibitem{villani2009optimal}
Villani, C.,
 \textit{Optimal transport: old and new},
  Vol. {338}, (Springer-Verlag, Basel, 2009).

\bibitem{Xu}
Xu, K., How has the literature on Gini's index evolved in the past 80 years? \emph{Dalhousie University, Economics Working Paper} (2003).  

\end{thebibliography}
\end{document}